\documentclass[preprint,12pt]{elsarticle}

\usepackage{amssymb}

\RequirePackage[OT1]{fontenc}
\RequirePackage{amsthm,amsmath, xcolor}
\RequirePackage[numbers]{natbib}
\RequirePackage[colorlinks,citecolor=blue,urlcolor=blue]{hyperref}
\usepackage{amscd,amsfonts,amssymb,latexsym,array,hhline,graphicx}
\usepackage{float}

\numberwithin{equation}{section}
\theoremstyle{plain}

\usepackage{amscd,amsfonts,amssymb,amsmath,latexsym,array,hhline,xcolor,graphicx}
\usepackage{float}

\newcommand\cE{{\cal E}}

\newcommand\cF{{\cal F}}

\newcommand\cJ{{\cal J}}

\newcommand\cL{{\cal L}}

\newcommand\cN{{\cal N}}

\newcommand\cD{{\cal D}}

\newcommand\cT{{\cal T}}

\newcommand\e{{\varepsilon}}

\def\E{{\bf E}}
\def\P{{\bf P}}

\def\bbr{{\mathbb R}}

\def\text#1{\hbox{#1}}

\def\E{{\bf E}}
\def\P{{\bf P}}

\def\R{{\bf R}}

\def\build #1_#2{\mathrel{\mathop{\kern 0pt #1}\limits_{#2}}} 
\newcommand{\zs}[1]{{\mathchoice{#1}{#1}{\lower.25ex\hbox{$\scriptstyle#1$}}
{\lower0.25ex\hbox{$\scriptscriptstyle#1$}}}}

\numberwithin{equation}{section}

\def\bbr{{\mathbb R}}

\def\bbr{{\mathbb R}}
\newcommand\fdem{$\Box$}
\newcommand\beq{\begin{equation}}
\newcommand\eeq{\end{equation}}
\newcommand\bea{\begin{eqnarray}}
\newcommand\eea{\end{eqnarray}}
\newcommand\bean{\begin{eqnarray*}}
\newcommand\eean{\end{eqnarray*}}

\newtheorem{theo}{Theorem}[section]
\newtheorem{prop}[theo]{Proposition}
\newtheorem{lemm}[theo]{Lemma}

\newtheorem{coro}[theo]{Corollary}

\begin{document}

\begin{frontmatter}

\title{Ruin Probabilities with Investments:
 Smoothness, IDE and ODE, Asymptotic Behavior}

\author[A1]{Yuri KABANOV } 
\author[A2]{Nikita  PUKHLYAKOV}
\address[A1]{Universit\'e Bourgogne Franche-Comt\'e, Laboratoire de Math\'ematiques,\\  16 Route de Gray,  25030 Besan\c{c}on cedex, France, and\\ 
Lomonosov Moscow State University, Russia\\
 Email: youri.kabanov@univ-fcomte.fr}

\address[A2]{Lomonosov Moscow State University, Russia\\
 Email: nikitapuhliakov@gmail.com}

\begin{abstract}
The study deals with the ruin problem when an insurance company 
having two business branches, life insurance and non-life insurance,   invests  its reserve into a risky asset with the price dynamics given by a geometric Brownian motion. 
We prove a result on smoothness of the ruin probability as  a function of the  initial capital and obtain for it an
integro-differential equation understood in the classical sense.   
For  the case of exponentially distributed jumps we show that the survival 
probability is a solution of an ordinary differential equation of the 4th order. Asymptotic analysis of the latter leads to the conclusion that  the ruin probability decays  to zero in the same way as in the already studied cases of models with one-side jumps.    
\end{abstract}



\end{frontmatter}

  \section{Introduction}

In the classical Lundberg--Cram\'er models of collective risk theory   insurance companies keep their reserve in cash (or in a bank account, typically, paying zero interest rate). In recent, more  realistic models, it is assumed that the capital reserves may vary not only due to the business activity but also due to a stochastic interest rate. In other words, an insurance company may invest all or just a part of its capital reserve into risky assets. These models lead to an important conclusion that the financial risk contributes enormously in the asymptotical behavior of the ruin probabilities: even under the Cram\'er condition, they are not exponentially decaying when the initial capital grows to infinity and  the ruin  always happens if the volatility of the risky asset is large with respect to the instantaneous interest rate. Moreover; they allows to quantify which proportion of risky investments may lead to the imminent ruin.   

Due to their  practical importance, ruin problems with investments became a vast and quickly growing chapter of the collective risk theory studying numerous models with various level of generality.  The ruin problem can be treated, at least, in two different  ways: using the techniques of integro-differential equations, see, e.g., \cite{Fr}, \cite{FrKP}, \cite{KP}, \cite{Belkina}, \cite{BelkinaKK},  or results from the  implicit renewal theory, see  
\cite{KP2019} and references wherein.  The first approach, allowing not only to obtain asymptotic but also  calculate ruin probabilities for given  values of the capital reserves, has some interesting mathematical questions.       

Our paper is a complement to the papers \cite{FrKP} and \cite{KP}, that extend, respectively,  the Lundberg--Cram\'er models for non-life insurance and 
life-insurance to the case where the capitals of insurance companies are invested into a risky asset with the price dynamics given by a geometric Brownian motion. In both,  the business activity  is given by  a compound Poisson process either with negative jumps  and positive drift (non-life insurance), or  
with positive jumps  and negative drift
 (life insurance). 
 Technically speaking,  these two models are quite different: in the first case the 
 downcrossing of zero may happen only at an instant of jump (thus, the model can be reduced to a discrete-time one) while in the second case the downcrossing happens in a continuous way and the reduction to a discrete-time model is not possible. Of course, in the classical setting if  jumps are exponentially distributed, a model with upward jumps can be transformed into a model with up downward jumps 
 and vice versa.  In models with investment the duality arguments does not work and one needs to treat them separately.   
   In the mentioned papers it was shown 
 that in the Lundberg-type model, that is,  with exponentially distributed jumps,  
 the ruin probability decreases with the rate $C u^{-\beta} $, $C>0$, when $\beta:=2a/\sigma^2-1>0$ and the ruin happens with probability one when $\beta\ge 1$ (the result for  $\beta=1$ was established in  \cite{PZ} and \cite{PerErr}). 
  
 Here we consider  a model of a company combining both types of business activities, that is, for which the corresponding compound Poisson process  has positive and negative jumps, \cite{Sax}.   
We prove,  under some minor assumptions  on the distributions of  jumps, that  the ruin probability as a function of the initial value is smooth  and satisfies an integro-differential equation (IDE). 
For a more specific case of exponential distributions we show that the ruin probability is a solution of the 4th order ordinary differential equation. Asymptotic analysis of this ODE leads us to our main result (Theorem \ref{main}) which looks exactly  as those 
of \cite{FrKP} and \cite{KP}.  Though the arguments follow the same general line as in  \cite{KP}, they are different in many aspects. Several seemingly  new results are obtained. 
In particular, for the model of Sparre Andersen  with investments, where the interarrival times form a renewal process,     
we derive a sufficient condition for the ruin with probability one (Theorem \ref{ruin}) and also  a lower asymptotic bound (Proposition \ref{Pr.sec:Lb.1}). For the model where the business activity of the company  can be represented as the difference of two compound Poisson processes with exponentially distributed jumps we prove a result on smoothness of the ruin probabilities.  


\section{The model}
\label{model}
Let  $(\Omega,\cF,{\bf F}=(\cF_t)_{t\ge 0},\P)$ be a stochastic basis, that is  a filtered probability space,  where we are given   a Wiener process $W=(W_t)_{t\ge 0}$ and an independent 
compound Poisson process $P=(P_t)_{t\ge 0}$  with drift $c$ and  successive jump instants $T_n$. 
 We denote $p_P(dt,dx)$ the jump measure of the latter with  its mean measure $\Pi_P(dx)dt$ where $\Pi_P(\R)<\infty$.  We assume that  $\Pi_P(\R_+)>0$ and $\Pi_P(\R_-)>0$ (some comments for  the cases where  $\Pi_P$ charges only of the half-axes will be also given). 

We consider  the process $X=X^{u}$ which  is the solution of   non-homogeneous linear stochastic 
equation 
\beq
 \label{risk}
 X_t=u+  \int_0^t  X_s  dR_s +P_t,
 \eeq
where $R_t=at+\sigma W_t$ is the relative price  process of the risky asset and $P$ is a compound Poisson process with drift $c\in \R$ representing the business activity of the insurance company, 
$u>0$ is the initial capital at time zero; we assume that  $\sigma^2>0$. 
The process  $X$  can be written in ``differential" form as $dX_t= X_t  dR_t +dP_t$, $X_0=u$, where  
\beq
 \label{r}
dP_t=cdt  + \int x  {p_P}(dt,dx), \quad P_t=0. 
\eeq 

In the actuarial context $X=X^u$ represents  the dynamics of the reserve  of an insurance company combining life and non-life insurance business and investing into a stock with the price given by  a geometric Brownian motion 
$$
S_t=S_0e^{(a-\sigma^2/2)t+\sigma W_t} 
$$
solving the linear stochastic differential equation $dS_t=S_tdR_t$ with initial value $S_0$. Changing in the need the money unit, we assume with our lost of generality that $S_0=1$. In this case we have by the product formula that 
\beq
\label{X}
X_t=S_t \Bigg (u +\int_0^t S_s^{-1}dP_s\Bigg). 
\eeq

In the classical collective risk theory the process $P$ are usually represented in the form 
\beq
\label{Pt}
P_t=ct+\sum_{i=1}^{N_t}\xi_i,  
\eeq
where $N_t:={p_P([0,t]\times \R)}$ is a Poisson process with intensity $\alpha=\Pi_P(\R)$
independent on the random variables
 $\xi_n=\Delta P_{T_n}$ where $T_n$ are successive instants of jumps of $N$; the i.i.d. random variables 
 $\xi_n$ has the probability distribution $F(dx)=\Pi_P(dx)/\Pi_P(\R)$.  Alternatively,  can represent $P$ using independent 
 compound Poisson processes given by the sums in the following representation.  
 $$
P_t=ct+\sum_{i=1}^{N^2_t}\xi^2_i-\sum_{i=1}^{N^1_t}\xi^1_i.  
$$
Here the  Poisson processes $N^1$ and  $N^2$ with intensities $\alpha_1=\Pi_P(\R_-)$  
and $\alpha_2=\Pi_P(\R_+)$, count, respectively, downward and upward jumps of $P$ and have successive jump instants $T^1_n$, $T^2_n$.  The  corresponding jump sizes $\xi^1_n$, $\xi^2_n$ are positive and have the  distribution functions   
$F_1(x):=\Pi_P(]-x,0] )/\alpha_1$,  $F_2(x):=\Pi_P(]0,x] )/\alpha_2$ for  $x>0$.

 Let $\tau^u:=\inf \{t: X^u_t\le 0\}$ (the instant of ruin),
 $\Psi (u):=P(\tau^u<\infty)$ (the ruin probability),
 and $\Phi(u):=1-\Psi (u)$
 (the survival probability).

\smallskip
Note that the same formula  (\ref{Pt}) $N$ can be an independent renewal process, that is, the counting process in which   
the lengths of the interarrival intervals $T_n-T_{n-1}$ form an i.i.d. sequence. 
In the collective risk theory this corresponds to the Sparre Andersen model.

\smallskip

The main aims of the present work: 

1) to get a result on smoothness in $u$ of the ruin probability $\Psi$ justifying  that  $\Psi$
solves  IDE in the classical sense; 

2) to deduce from the latter, in the special case of exponentially distributed jumps, an ODE  and use it to obtain the  exact asymptotic of the ruin probability as $u\to \infty$. 
\smallskip

Of course, in the model with only upward jumps (i.e., with $\Pi_P(]\infty,0[)= 0$) if $c\ge 0$  there is no ruin; the nontrivial  case 
where $c<0$ is studied in details in \cite{KP} and \cite{PZ}. The paper \cite{FrKP} deals with the model with downward jumps, i.e. with the non-life insurance, but the question of smoothness is not discussed in and, seemingly, needs to be revisited (in \cite{WangWu} the smoothness is established under restrictions on parameters).

\smallskip 
\noindent
{\bf Notations.} Throughout the paper we shall use the following abbreviations: 
$$
\qquad \beta:=
 {2a}/{\sigma^2}-1, \qquad \kappa:=a-\sigma^2/2=(1/2)\sigma^2\beta, \qquad \eta_t:=\ln S_t=\kappa t+\sigma W_t.   
$$  

The following simple result holds for any Sparre Andersen model with investments. 

\begin{lemm} Suppose that there is $\beta'\in ]0,\beta\wedge 1[$ such that $E(\xi^-_1)^{\beta'}<1$. Then $\Psi(u)\to 0$ as $u\to \infty$. 
 \end{lemm}
 {\sl Proof.} Let 
 $\tilde \Psi(u)$ be the ruin probability for the reserve process  $\tilde X^u$ corresponding to the model where the business of the company is given by   
$$
\tilde P_t=-|c|t-\sum_{i=1}^{N_t}\xi^-_i.  
$$
Then $\Psi(u)\le \tilde \Psi(u)\le \P(Z_{\infty}>u)$ where 
$$
Z_{\infty}:= -\int_0^{\infty} S_s^{-1}d\tilde P_s=|c|\int_0^{\infty}  e^{-\kappa s - \sigma W_s}ds + \sum_{n=1}^\infty e^{-\kappa T_n-\sigma W_{T_n}}\xi_k^-. 
$$
Then $\kappa s + \sigma W_s=(1/2)\sigma^2s(\beta+ (2/\sigma)W_s/s)$. By the law of large numbers for almost all $\omega$ there is $s_0(\omega)$ such that  $\beta+ (2/\sigma)W_s/s>\beta/2$ when $s\ge s_0(\omega)$.  Thus,  the integral in the right-hand side above is a finite random variable. Also 
$$
e^{-\kappa T_n-\sigma W_{T_n}}=\prod_{k=1}^n \zeta_k,
$$
where $\zeta_k:=e^{-\kappa (T_{k} -T_{k-1})-\sigma (W_{T_{k}}-W_{T_{k-1}})}$  form an i.i.d. sequence.  Note that $\E\zeta_1^\beta=1$. Hence, $\E\zeta_1^{\beta'} <1$ and  
$$
\E\sum_{n=1}^\infty(\xi^-_n)^{\beta'} \prod_{k=1}^n \zeta^{\beta'}_k <\infty.
$$
That is, $\sum_{n=1}^\infty(\xi^-_n \prod_{k=1}^n \zeta_k)^{\beta'} <\infty$ (a.s.). But then  $\sum_{n=1}^\infty\xi^-_n\prod_{k=1}^n \zeta_k <\infty$ (a.s.). \fdem

%
%
 \begin{theo}
\label{main}
  Let $F_1(x)=1-e^{-x/\mu_1}$ and  $F_2(x)=1-e^{-x/\mu_2}$ for $x>0$.  Assume that $\sigma>0$ and $P$ is not an increasing process.

 (i) If
  $\beta>0$, then for some $K>0$
 \beq
  \Psi (u)= Ku^{-\beta}(1+o(1)), \quad u\to \infty.
  \eeq

(ii) If $\beta\le 0$, then $\Psi (u)=1$ for all $u>0$.
 \end{theo}

\section{Large volatility case: ruin with probability one}
\label{secruin}
The result below gives a sufficient condition on the ruin with probability one for Sparre Andersen models with risky investments. It implies the statement $(ii)$ of Theorem \ref{main}. 

\begin{theo} 
\label{ruin}
Suppose that $P$ is a non-increasing compound renewal process with drift (i.e., $c<0$ or $\P(\xi_1<0)>0$) such that    
$\E|\xi|^\e<\infty$ and $\E e^{\e T_1}<\infty$ for some $\e>0$.  
If $\beta\le 0$,  then  $\Psi(u)=1$ for any $u>0$. 
\end{theo}

{\sl Remark.}  In this formulation one can replace the assumption $\E|\xi|^\e<\infty$ by  a formally weaker $\E(\xi^+)^\e<\infty$.

  As in \cite{KP} the arguments are based on the   ergodic property of the  discrete-time autoregressive 
process $(\tilde X_n^u)_{n\ge 1}$ with random coefficients  which is defined   recursively by the relations 
\beq
\label{recur}
\tilde X_n^u=A_n\tilde X_{n-1}^u +B_n, \qquad n\ge 1, \quad \tilde X_0^u=u,  
\eeq
where $(A_n,B_n)_{n\ge 1}$ is an  i.i.d.  sequence  in 
$\bbr^2$. For the following result see \cite{PZ}, Prop. 7.1.  
\begin{lemm} 
Suppose that   $\E |A_1|^\delta<1$ and 
$\E |B_1|^\delta<\infty$ for some $\delta\in ]0,1[$. Then for any $u\in \bbr$ the sequence $\tilde X_n^u$ converges in $L^\delta$ (hence, in probability) to the random variable 
$$
\tilde X_\infty^0=
\sum_{k=1}^\infty B_k\prod_{j=1}^{k-1}A_j
$$  
and for any bounded uniformly continuous function $f$
\beq\label{ergodic}
\frac 1n\sum_{k=1}^n f(\tilde X_k^u)\to \E f(\tilde X_\infty^0) \quad \hbox{in probability as } n\to \infty . \eeq  
\end{lemm}


\begin{coro} 
\label{coro}
Suppose that   $\E |A_1|^\delta<1$ and $\E |B_1|^\delta<\infty$ for some $\delta\in ]0,1[$. 

$(i)$ If $\P(\tilde X_\infty^0<0)>0$, then  $\inf_{k\ge 1}\tilde X_k^u<0$. 

$(ii)$ If $A_1>0$   and  $B_1/A_1$ is unbounded from below, then   $\inf_{n\ge 1}\tilde X_k^u<0$. 
\end{coro}

\noindent
{\sl Proof.} $(i)$  Let $f(x):=-I_{\{x<-1\}}+x I_{\{-1\le x<0\}}$. Then $\E f(\tilde X_\infty^0)<0$ and (\ref{ergodic}) may hold only if  $\inf_{k\ge 1}\tilde X_k^u<0$. 

$(ii)$ Put  $\tilde X_\infty^{0,1}:=\sum_{n= 2}^\infty B_n\prod_{j=2}^{n-1}A_j$. Then 
$$
\tilde X_\infty^0=B_1+A_1\, \tilde X_\infty^{0,1}=A_1(\tilde X_\infty^{0,1}+B_1/A_1). 
$$ 
Since $B_1/A_1$ and $\tilde X_\infty^{0,1}$ are independent random variables and  $B_1/A_1$ is unbounded from below, so is the sum $\tilde X_\infty^{0,1}+B_1/A_1$. 
It follows that the probability  $\P(\tilde X_\infty^0<0)>0$ and we can use $(i)$.   
\fdem

\medskip
\noindent
{\sl Proof of Theorem \ref{ruin}.}
Let $\tilde X_n=\tilde X^u_n:=X^u_{T_n}$ where $X^u$ is the process given by 
the formulae (\ref{X}) and (\ref{Pt}) assuming that $N$ is a renewal process. 
In this case 
$$
X^u_{T_n}=e^{\eta_{T_n}}u+\sum_{k=1}^ne^{\eta_{T_n}-\eta_{T_k}}\left(c\,\int^{T_{k}}_{T_{k-1}}\, e^{\eta_{T_{k}}-\eta_{v}}d v+\xi_{k}
\right)
$$
and 
\beq
\label{MQ}
A_{k}:=e^{\eta_{T_{k}}-\eta_{T_{k-1}}},  
\qquad
B_{k}:=\xi_{k}+c\,\int^{T_{k}}_{T_{k-1}}\, e^{\eta_{T_{k}}-\eta_{v}}d v.
\eeq
Since the Wiener process $W$ and  $P$  (a compound  renewal process with drift) are independent, 
$(A_k,B_k)_{k\ge 1}$ is an i.i.d. sequence.  

According to the formulae (\ref{X}) and (\ref{Pt})
\beq
\label{tildes}
\tilde X_{n}=\cE_n u
+\sum^{n}_{k=1}B_{k}\frac {\cE_n}{\cE_k}, \qquad  \cE_n:=\prod^{n}_{j=1}\,A_{j}.
\eeq
Clearly, $\tilde X_{n}=A_n\tilde X_{n-1}+B_n$, that is, $\tilde X_{n}$ satisfies 
(\ref{recur}).  

\begin{lemm}
\label{unbb}
Suppose that $c<0$ or $\P(\xi_1<0)>0$. Then 
the ratio 
$$
\frac{B_1}{A_1}=e^{-\eta_{T_1}}\xi_1+c\int_0^{T_1}e^{-\eta_s}ds
$$ 
is unbounded from below on the (non-null) set $\{W_{T_1}<0\}$. 
\end{lemm}
{\sl Proof}. Recall that the conditional law of the the Wiener process $W$ on $[0,t]$ given $W_t=x$ is the same as the Brownian bridge  $[0,t]$ ending at $x$, i.e. 
coinciding with the law of the process $(W_s-(s/t)W_t+(s/t)x)_{s\le t}$.  

Fix $N>0$. Let $c<0$. 
Then 
\beq
\label{fond1}
\P(B_1/A_1\le -N|(T_1,W_{T_1})=(t,x))=\P\left(e^{-\kappa t-\sigma x}\xi_1+ c\zeta^x_t\le -N\right)
\eeq
where the random variable 
$$
\zeta^x_t:=\int_0^t e^{-\kappa s-\sigma (W_s-(s/t)W_t+(s/t)x)}ds
$$ 
is independent on $W$. The law of $\zeta^x_t$ charges every open interval in 
$\bbr_+$. Thus, the right-hand side  of (\ref{fond1}) is strictly positive for every $(t,x)\in ]0,\infty[\times \bbr$. 

Let $\P(\xi_1<0)>0$. Take $\e>0$ and $T>0$ such that the probabilities $\P(\xi_1\le -\e)$ and $P(T_1\le T)$ are strictly positive. Then 
$$
\P(B_1/A_1\le -N, \xi\le -\e |(T_1,W_{T_1})=(t,x))\ge \P\big(-\e e^{-\kappa T-\sigma x} + c\zeta^x_t\le -N\big)>0
$$
for all sufficiently small $x$, namely, satisfying the inequality $\e e^{-\kappa T-\sigma x}>N$.   Since the law  of $(T_1,W_{T_1})$  charges the set $[0,T]\times ]-\infty,y]$ whatever is $y\in \bbr$, the result follows. \fdem 

\begin{lemm} 
\label{momb}
Let $\beta\le 0$. If $\E|\xi|^\e<\infty$ and $\E e^{\e T_1}<\infty$ for some $\e>0$, then $\E |B_1|^\delta<\infty$ for some $\delta\in ]0,1[$. 
\end{lemm}
{\sl Proof.}  Note that  $\eta_{T_1}-\eta_s\le \sigma (W_{T_1}-W_s)$. For any  $\delta\in ]0,1[$  
$$
\E|B_1|^{\delta}\le 
\E|\xi_1|^{\delta}+|c|^\delta\E\left(\int_0^{T_1}e^{\eta_{T_1}-\eta_s}ds\right)^{\delta}\le  
\E|\xi_1|^{\delta}+|c|^\delta \E T_1^\delta\sup_{s\le T_1} e^{\delta \sigma W_s}. 
$$
 Let $m(dt)$ denote the distribution of $T_1$. Substituting the density of distribution of the running maximum of the Wiener process we get that  
\bean
\E T_1^\delta \sup_{s\le T_1} e^{\delta \sigma W_s}&=&\int_0^\infty t^\delta \int_0^t \sqrt{\frac 2{\pi t}}e^{\delta \sigma x}
e^{-x^2/(2t)}dxm(dt)\\
&=&
\int_0^\infty t^\delta  e^{(1/2)\delta^2\sigma^2 t}\int_0^{\sqrt t}\sqrt{\frac 2{\pi}}e^{-(y-\delta\sigma\sqrt t)^2/2} 
dym(dt)\\
&\le & 2\E T_1^\delta e^{(1/2)\delta^2\sigma^2 T_1}. 
\eean
Our assumptions imply that  $\E|B_1|^{\delta}<\infty$ for all sufficiently small $\delta>0$.  \fdem

\smallskip
In the particular case where $\beta<0$ we can get the result by direct reference 
to Corollary \ref{coro}$(ii)$ because  
$$
\E A_1^{-\beta}= \E e^{-\beta\eta_{T_1}}=\int_0^\infty \E e^{-\beta\sigma W_t-(1/2)(\beta\sigma)^2t}m(dt)=1 
$$ 
and, therefore,   $\E A_1^{\delta}<1$ for any $\delta\in ]0,-\beta[$. 
 
%

To cover the general case, we consider a suitably chosen random subsequence of 
$\tilde X_n$ which satisfies a linear difference equation with needed properties. 

Let $\widehat X_n=\widehat X^u_n:=\tilde X_{\theta^n}$ where 
$\theta_n:= \inf\{k> \theta_{n-1}\colon  \cE_k<\cE_{\theta_{n-1}}\}$. Thus,  $\theta_n$ are ladder times for the random walk $M_k=\ln \cE_k$. If $\beta=0$, then 
$M_1=\sigma W_{T_1}$ and $ \E M_1=0$, $ \E M_1^2=\E T_1<\infty$. 
 Therefore,   
 there is a finite constant $C$ such that 
\begin{equation}\label{sec:Ras.5}
\P(\theta_{1} > n) \le Cn^{-1/2},
\end{equation}
see Theorem 1a in Ch. XII.7 of Feller's book \cite{Fel} and the remark before it. In the general case, $M_1=(1/2)\beta \sigma^2 T_1+\sigma W_{T_1}$ and the above bound holds also when $\beta<0$.   

In particular, $\theta_n<\infty $ and the differences $\theta_n-\theta_{n-1}$ form  a sequence of finite  independent random variables distributed as $\theta_1$. The discrete-time process 
$$
\widehat X^u_n=\cE_{\theta_n}u+\sum^{\theta_n}_{k=1}B_{k}\frac {\cE_{\theta_n}}{\cE_k}
$$ 
solves the linear equation 
$$
\widehat X_n^u=\widehat A_n\widehat X_{n-1}^u +\widehat B_n, \qquad n\ge 1, \quad \widehat X_0^u=u,  
$$
where
$$
\widehat  A_n:=\frac {\cE_{\theta_n}}{\cE_{\theta_{n-1}}}, \qquad \widehat B_n:= 
\sum_{k=\theta_{n-1}+1}^{\theta_n} B_k\frac {\cE_{\theta_n}}{\cE_k}. 
$$


By construction, $ \widehat A_1<1$ and  
$$
\vert \widehat B_{1} \vert\le \sum_{k=\theta_{n-1}+1}^{\theta_n}\vert B_k\vert \frac {\cE_{\theta_n}}{\cE_k}
\le  \sum^{\theta_{1}}_{j=1} 
\,\vert {B}_{j}\vert.
$$
According to Lemma \ref{momb}   $\E\vert B_{1}\vert^{\delta}<\infty$ 
 for some $\delta\in ]0,1[$.
Taking  $r\in ]0,\delta /5[$ and defining the sequence $l_{n}:=[n^{4r}]$, we have, using the Chebyshev inequality and  (\ref{sec:Ras.5}), 
that 
\begin{align*}
\E\vert \widehat B_{1} \vert^{r}&\le 1+r\sum_{n\ge 1}\,{n^{r-1}}\P\Bigg(  \sum^{\theta_{1}}_{j=1}
\,\vert B_{j}\vert
>n
\Bigg)\\[2mm]
&\le  1+r\sum_{n\ge 1}\,{n^{r-1}}\,\P\Bigg( \sum^{l_{n}}_{j=1}
\,\vert B_{j}\vert>n\Bigg)+ r\sum_{n\ge 1}\,{n^{r-1}}
\P(  \theta_{1}>l_{n})
\\[2mm]
&\le 
 1+r\E\vert Q_{1}\vert^{\delta}
 \sum_{n\ge 1}\,l_{n}n^{r-1-\delta} +rC\sum_{n\ge 1}\,n^{r-1}l_{n}^{-1/2}<\infty.
\end{align*}
 
To apply Corollary \ref{coro}$(ii)$ it remains to check that the random variable  $\widehat B_1/\widehat A_1$ is unbounded from below.  
But Lemma \ref{unbb} asserts that this ratio coinciding with $B_1/A_1$ on the set $\{W_{T_1}<0\}$ of strictly positive probability is unbounded from below on this set.  \fdem

 \section{Lower asymptotic bound}
\label{sec:Lb}
The next result we are needed at our asymptotic analysis  indicates  that the ruin probability decreases at infinity not faster than a certain power function.  The proof given below covers also the more general case  where $P$
is a compound renewal process with drift given by the representation (\ref{Pt}) where $N$ is a counting renewal process. 

 \begin{prop}
  \label{Pr.sec:Lb.1}
 Suppose that $c<0$ or the random variable $\xi_1$ is unbounded from below. Then there exists $\beta_{*}>0$ such that
 \beq
 \label{lab}
 \liminf_{u\to\infty}\,u^{\beta_{*}}\,\Psi(u)\,>\,0\,.
 \eeq
 \end{prop}
 \noindent {\sl Proof.}  
 Let $\tilde X_n=\tilde X^u_n:=X^u_{T_n}$ and let  $\theta^u:=\inf\{n:\ \tilde X^ u_n\le 0\}$.  If $c<0$ then the ruin may happen between jump times but in all cases, 
 $$
 \Psi(u):= \P(\tau^u<\infty) \ge \P(\theta^u<\infty).
 $$ 

Recall that for $(\tilde X_n)$ we have the formulae (\ref{tildes}) and (\ref{recur}) with 
$(A_k,B_k)$ defined by (\ref{MQ}). 


For  reals  $\varrho\in ]0,1[$ and $b>1/({\varrho^{2}(1-\varrho)})$ we
define the sets  
\begin{equation}\label{sec:Lb.2}
\Gamma_{k}:=\{A_{k}\le \varrho\}\cap\{B_{k}\le \varrho^{-1}\}\,,
\quad
D_{k}:=\{A_{k}\le \varrho^{-1}\}\cap\{B_{k}\le -b\}.
\end{equation}
Note that $\P (\Gamma_{k})=\P (\Gamma_{1})$ and $\P (D_{k})=\P (D_{1})$ for all $k$. 

\begin{lemm}
If there are $\varrho$ and $b$ such that $\P (\Gamma_{1})>0$ and $\P (D_{1})>0$, then (\ref{lab}) holds. 
\end{lemm} 
\noindent 
{\sl Proof.} 
Using (\ref{tildes}) we easily get that  on the set  $\cap^{n}_{k=1}\,\Gamma_{k}$
$$
\tilde X_{n}\le  u \varrho^{n}+\frac{1}{\varrho(1-\varrho)}.
$$
From  the representation  $\tilde X_{n+1}=A_{n+1} \tilde X_{n} + A_{n+1}$ and the above bound we infer that  on the set 
$\big(\cap^{n}_{k=1}\,\Gamma_{k}\big)\cap D_{n+1}$
$$
\tilde X_{n+1}\le u \varrho^{n-1}+\frac{1}{\varrho^{2}(1-\varrho)}
-b= u \varrho^{n-1} -b_1.
$$
where $b_1:=b-1/({\varrho^{2}(1-\varrho)})>0$. 

Let  $u>b^1$ and let $n=n(u):=2+[{(1/\ln\varrho)}{\ln(b_{1}/u)}] $ where $[x]$ means the integer part, 
$x-1<[x]\le x$. 
 Then  $$
u \varrho^{n-1}=  u e^{(n-1)\ln \varrho}< ue^{\ln (b_1/u)} < b_1
$$
 and, therefore,  
%
%
$$
\P(\theta^u<\infty)\ge \P\left( \cap^{n}_{k=1}\,\Gamma_{k}\cap D_{n+1}\right)
=\left(\P(\Gamma_{1})\right)^{n}\P(D_{1}).
$$
Take $\beta_*:=\frac{\ln\P(\Gamma_{1})}{\ln\varrho}$. Then 
$$
 u^{\beta_*} \P(\theta^u<\infty)\ge  e^{\frac {\ln\P(\Gamma_{1})}{\ln\varrho}\ln u+n\ln \P(\Gamma_{1})}\P(D_{1})
 \ge e^{\big(2+\frac {\ln b_1}{\ln \varrho}\big)\ln \P(\Gamma_{1})}\P(D_{1}). 
$$
Since $\P(\tau^u<\infty)\ge \P(\theta^u<\infty)$, the lemma is proven. \fdem 
\smallskip

It remains to show that our assumptions ensure that for any $ \varrho, b>0$  
the probability $\P (A_{1}\le \varrho, \ B_{1}\le -b)$ is strictly positive.  
We use  the fact that 
the conditional distribution of the Wiener process $(W_s)_{s\le t}$ given $W_t=x$ is the (unconditional) distribution  of the Brownian bridge $(B^x_s)_{s\le T}$ which is the same the distribution of the process  $(W_s-(s/t)W_t)_{s\le t}$.  
Take $v_1, v_2>0$ and  $x_0$  such that $\P(T_1\in [v_1, v_2])>0$ and  $|\kappa |v_2 + \sigma x_0 \le \ln \varrho$. Then for any $(t,x)\in \Delta:= [v_1, v_2]\times ]-\infty,x_0]$ we have that  
\beq
\label{Ups}
\P (A_{1}\le \varrho, \ B_{1}\le -b\;\vert\; (T_1,W_{T_1})=(t,x))=\P(\xi_1+ c\Upsilon ^x_{t}\le -b)
\eeq
where 
$$
\Upsilon^x_t:=\int_0^te^{\kappa (t-s) + \sigma x  - \sigma(W_s-(s/t)W_t)}ds.
 $$
Since the distribution of Wiener process has a full support, i.e.   charges any open set in the space $C_0([0,T])$ of the trajectories, the random variable  $\Upsilon^x_t$  is unbounded from above (whatever are 
$x$ and $t>0$). Any of our assumptions implies that $\xi_1+ c\Upsilon ^x_{t}$ is unbounded from below, that is the probability in the right-hand side of (\ref{Ups}) is strictly positive. Integrating (\ref{Ups}) over $\Delta$ with respect to the distribution of $(T_1,W_{T_1})$ charging $\Delta$ we get the result.  \fdem


\section{Regularity of the ruin probability}
\label{sec:Reg}

Following tradition (seemingly, justified by  notational convenience) we shall work with the survival probability $\Phi=1-\Psi$ which has the same regularity property and satisfies the same 	equations. 						
 \begin{theo}
\label{Th.sec:Reg.1}
Suppose that  $F_k(dx)=f_k(x)dx$  where the densities $f_k$ are two times  differentiable on $]0,\infty[$ and   $f'_k,f''_k\in L^1(\bbr_{+})$, $k=1,2$.
 Then $\Phi $ is two times continuously differentiable on $]0,\infty[$ and $\Phi'$, $\Phi''$ are bounded. 
 \end{theo}
 
%
\noindent 
{\sl Proof.} Define   the  continuous  process
\begin{equation}
\label{sec:Reg.0}
Y^{u}_{t}:=S_t\left (u+c\int_{[0,t]} S^{-1}_s ds\right)
\end{equation}
coinciding with $X^u$ on $[0,T_1[$ 
 and introduce  the stopping time 
 $$
 \theta^{u}:=\inf\{t\ge 0\colon \ Y^{u}_{t}\le 0\}.
 $$
 By virtue of  the strong Markov property of $X ^u$
\begin{equation}
\label{sec:Reg.0-1}
\Phi(u)=\E\Phi(X^u_{\theta^{u}\wedge T_1})=\E\Phi(Y^u_{\theta^{u}\wedge T_1}+\Delta X^u_{\theta^{u}\wedge T_1}).
 \end{equation} 
Due to independence of $W$ and  the  Poisson processes $N^1$, $N^2$,  the values of $\theta^u(\omega) $,
  $T^1_1(\omega)$ and $T^2_1(\omega)$ are all different for almost all  $\omega$. 
Since $\Phi(X^u_{\theta^{u}})I_{\{\theta^{u}< T_1\}}=0$ we have the representation 
$\Phi=K^\downarrow +K^\uparrow $ where 
$$
K^\downarrow (u):=\E I_{\{Y^u_{T^1_1}>0\}}
I_{\{T_1^1<T_1^2\}}\Phi(Y^u_{T^1_1}-\xi^1_1),  \quad  K^\uparrow (u):=\E I_{\{Y^u_{T_1^ 2}>0\}}
I_{\{T_1^1>T_1^2\}}\Phi(Y^u_{T_1^2}+\xi^2_1). 
$$
The analysis of the smoothness of these two functions is similar, so we consider the first one.  
It is convenient to represent it as the sum $K^\downarrow =K_1^\downarrow +K_2^\downarrow $ where 
\bea 
\label{K1}
K_1^\downarrow (u)&:=&\int_{\R^3}  I_{\{0<s<t\wedge 2\}}\E G (Y^u_s,w)m(ds,dt)\,F_1(dw),\\
\label{K2}
K_2^\downarrow (u)&:=&\int_{\R^3}  I_{\{2<s<t\}}\E G(Y^u_s,w)m(ds,dt)\,F_1(dw), 
\eea
with
$G(y,w):=I_{\{y> 0\}}\Phi\left(y-w\right)=I_{\{y-w> 0\}}\Phi\left(y-w\right)=\Phi\left(y-w\right)$ for $w\ge 0$ and  
$$
m(ds,dt):=\alpha_1\alpha_2e^{-(\alpha_1s+\alpha_2 t)}ds\,dt.
$$


\begin{lemm}
For any bounded measurable function $G(y,w)$ the function $K_2^\downarrow (u)$ defined  (\ref{K2}) belongs to $C^{\infty}(]0,\infty[)$ and has bounded derivatives of any order. 
\end{lemm}
%
%
\noindent
{\sl Proof.}
Using the representation 
$$
Y_s^u= e^{\eta_s-\eta_1}Y_1^u+\int_{[1,s]} e^{\eta_s-\eta_r} dr, \qquad s\ge 2,  
$$
and noticing that the random variable  $Y_1^u$ is independent of the process 
$(\eta_s-\eta_1)_{s\ge 1}$, we get that  
$$
\E (G(Y^u_s,w)|Y^u_1)=G(s,Y_1^u,s,w)
$$
where 
$$
G(s,y,w):=\E G\Bigg(e^{\eta_s-\eta_1}y+\int_{[1,s]} e^{\eta_s-\eta_r} dr,w\Bigg). 
$$
Substituting the expression for $Y_1^u$ we obtain the formula  
$$
\E G(Y^u_s,w)=\E G(s,Y^u_1,w)=\E G(s,e^{\kappa+\sigma W_1}(u +c R_1),w)
$$
with  
$$
R_1:=\int_{[0,1]}e^{-\kappa r -\sigma W_r}dr.  
$$
Recall that the conditional distribution of the process $(W_r)_{r\le 1}$ given $W_1=x$, i.e. the unconditional distribution of the Brownian bridge $B^x$ with $B^x_0=0$ and $B^x_1=x$, which  is the same as of the process 
$W_r+r(x-W_1)$, $r\le 1$. It follows that 
$$
\E G(Y^u_s,w)=\int_{\R} \E G\big(s,e^{\kappa +\sigma x}(u+ \zeta^x),w\big)\varphi_{0,1}(x)dx
$$
where the random variable 
$$
\zeta^x:=c\int_{[0,1]}e^{-\kappa r-\sigma (W_r+r(x-W_1))}dr. 
$$

The case  $c= 0$ is easy. By  the change of variable $z=\kappa +\sigma x+\ln u$ we get that 
$$
\E G(Y^u_s,w)=\frac 1\sigma \int_{\R} G(s,e^z,w)\varphi_{0,1}((z - \kappa- \ln u)/\sigma)dz. 
$$
The function $u\mapsto \varphi_{0,1}((z - \kappa- \ln u)/\sigma)$ belongs to $C^{\infty}(]0,\infty[)$
and its derivatives on any interval $[u_1,u_2]\subset ]0,\infty[$ are dominated by integrable functions. It follows that  the function 
$u\mapsto \E G(Y^u_s,w)$ belongs to $C^{\infty}(]0,\infty[)$  and its derivatives are locally bounded. 

\medskip
Let  $c\neq 0$. 
Lemma \ref{F2}  below asserts that  for every  $x$ the random variable $\zeta^x$ has a density $\rho (x,.)$  so that 
$$
\E G(Y^u_s,w)=\int_{\R^2} G\big(s,e^{\kappa +\sigma x}z,w\big)\rho(x,z-u)\varphi_{0,1}(x)dxdz.
$$
Since this density belongs to $C^{\infty}$ and its derivatives are of sub exponential growth in $x$ the function  $y\mapsto \E G(Y^u_s,w)$ also belongs to $C^{\infty}$ and has bounded derivatives.  
So, the same property has the function $u\mapsto K_2^\downarrow (u)$ and the lemma is proven.

\begin{lemm} [see \cite{KP}, Lemma 5.2]
\label{F2}
The random variable $\zeta^x$ has a density  $\rho(x,.)\in C^\infty$ such that  for any $n\ge 1$  
\beq
\label{bderivs}
\sup_{y\ge 0}\,
\left\vert \frac {\partial^n}{\partial y^n} \rho(x,y) \right\vert  \le C_n e^{C_n |x|}
\eeq
with some constant $C_n$ and  $(\partial^n /\partial y^n)\rho(x,0)=0$. 
\end{lemm}
The needed smoothness property of $K_1^\downarrow$ follows from the following lemma. 
\begin{lemm} Let $\xi>0$ be a random variable with a density $f$ which is two times  differentiable on $]0,\infty[$ and   $f',f''\in L^1(\bbr_{+})$.  Let   $G: \R\to [0,1]$ be a measurable function vanishing on $]-\infty,0]$ and let $h(y):= \E G(y-\xi)$. Then  the function $(s,u)\mapsto \E h(Y_s^u)$ has two continuous derivatives in $u$ bounded on $[0,t]\times ]0,\infty]$. 
\end{lemm}
\noindent
{\sl Proof.}
First, we observe that 
$$
h(y)= \int_\R G(y-x)f(x)dx =\int_\R G(z)f(y-z)dz
$$
and 
$$
h'(y)=\int_\R G(z)f'(y-z)dz. 
$$
It follows that  $|h'(y)|\le ||f'||_{L^1}$ and $|h''(y)|\le ||f''||_{L^1}$.

Using the representation
$$
Y^u_s=e^{\eta_s}u+\int_{[0,s]}e^{\eta_s-\eta_r}ds 
$$ 
and arguing in the same spirit as above but conditioning this time on the random variable $W_s\sim \cN(0,s)$ and considering the Brownian bridge on $[0,s]$  we obtain that 
$$
\E G(Y_s^u-\xi)=\frac 1{\sqrt s}\int_{\R}\E h(e^{\kappa s+\sigma x }(u+\zeta^{s,x}))\phi_{0,1}(x/\sqrt{s})dx
$$
where 
$$
\zeta^{s,x}:=c\int_{[0,s]}e^{- (\sigma r x/s+\kappa r+
\sigma (W_{r}-(r/s)W_{s})
}dr. 
$$
If $c=0$, then  
$$
\E G(Y_s^u-\xi)=\int_{\R}h(e^{\kappa s+\sigma {\sqrt s}x }u)\phi_{0,1}(x)dx
$$ 
and the needed property is obvious. 
%

It is easily seen that the random variable $\zeta^{s,x}$ has a $C^\infty$ density  (the same as of 
$\zeta^x$ but with the parameters $cs$, $\kappa s$, and $\sigma s^{1/2}$). 
Unfortunately,  derivatives of this density have non-integrable singularities at $s=0$.   By this reason  arguments used above do not work. 

Let $c\neq 0$. Then 
the smooth function  $x\to \zeta ^{s,x}$ is strictly decreasing  and maps $\bbr$ onto $\bbr_+$ (when $c>0$) or strictly increasing  and maps $\bbr$ onto $\bbr_-)$ (when $c<0$). Let denote 
 $z(s,.)$ its inverse which is a function decreasing from $+\infty$ to $-\infty$ (when $c>0$) and 
 increasing from $-\infty$ to $+\infty$ (when $c<0$). The  partial   derivative in $x$ is given by the formula
\beq
\label{zx}
z_x(s,x)=-\frac{s}{L(s,z(s,x))},
\eeq
where
\begin{equation}
\label{sec:Reg.7}
L(s,x)=c \sigma \int_{[0,s]}r e^{- (\sigma r x/s+\kappa r+
\sigma (W_{r}-(r/s)W_{s})
}dr.
\end{equation}
In both cases $z_{xx}(s,x)>0$ for $s>0$. 

Changing the variable, we obtain that $\E h(Y_s^u)=\E H(s,u)$ where 
\beq
\label{sec:Reg.9a}
H(s,u):=
\frac 1{\sqrt s}\int^{\infty}_{0}
h(e^{\kappa s+\sigma z(s,x) }(u+x)){\varphi_{0,1}\left({z(s,x)}/\sqrt{s}\right)}z_x(s,x)d x, \qquad c>0, 
\eeq
For $s\in ]0,2]$
\bean
\frac 1{\sqrt{s}}\int^{\infty}_{0}e^{\sigma z(s,x) }
\varphi_{0,1}\left({z(s,x)}/\sqrt{s}\right) z_x(s,x)d x&=&\int e^{\sigma \sqrt{s}z}\varphi_{0,1}(z)dz\\
&\le &
\int e^{\sigma \sqrt{2}z}\varphi_{0,1}(z)dz
\eean
where the last integral is a finite constant. This estimate and the boundedness of $h'$ legitimate the differentiation under the sign of the integrals.   In the case  $c>0$ we have  the formula 
$$
H_u(s,u)=\frac 1{\sqrt{s}}\int^{\infty}_{0}e^{\kappa s+\sigma z(s,x) }
h'(e^{\kappa s+\sigma z(s,x) }(u+x)){\varphi_{0,1}\left({z(s,x)}/\sqrt{s}\right)}z_x(s,x)d x. 
$$
and the bound $|H_u(s,u)|\le C$ for all $s\in ]0,2]$. 

Repeating the arguments we get that 
$$
H_{uu}(s,u)=\frac 1{\sqrt{s}}\int^{\infty}_{0}e^{2\kappa s+2\sigma z(s,x) }
h'(e^{\kappa s+\sigma z(s,x) }(u+x)){\varphi_{0,1}\left({z(s,x)}/\sqrt{s}\right)}z_x(s,x)d x. 
$$
$|H_u(s,u)|\le C$ for all $s\in ]0,2]$. 

Similar arguments are applied  in the case $c<0$.  
\fdem

Theorem \ref {Th.sec:Reg.1} is proven. \fdem 
\smallskip 

\noindent 
{\sl Remark.} 
Let $V: \bbr\to [0,1]$ be a measurable function, $V(x)=0$ for $x\le 0$. Put 
$$\Psi_V(u):=EV(X^u_{\tau^u})I_{\{\tau^u<\infty\}}).
$$ 
Then the statement of Theorem \ref{Th.sec:Reg.1} holds for $\Psi_V$ with the same proof. Indeed, 
the strong Markov property for $\Psi_V$ has the same form as for $\Phi$ in (\ref{sec:Reg.0-1}) that is $\Psi_V(u)=\E\Phi_V(X^u_{\theta^{u}\wedge T_1})$. 
Also Proposition \ref{propos1} below holds for $\Psi_V$. 
In the particular case where  $V(x)=1$ for all $x<0$, the function   $\Psi_V$ coincides with $\Psi$ on $]0,\infty[$ but  they are different on $\bbr_-$.

\section{The integro-differential equation for the survival probability}
\begin{prop}
\label{propos1}
Suppose that $\Psi\in C^2$. Then the function $\Phi$ on $]0,\infty[$
satisfies the following equation:
\beq
\label{int-diffG}
 \frac 12 \sigma^2u^2 \Phi''(u)+ (au+c)\Phi'(u) +\int (\Phi(u+y)-\Phi(u))\Pi_P(dy)=0.
\eeq
 \end{prop}
{\sl Proof.} 
For $h>0$ and $\epsilon>0$  small enough to ensure that $u\in ]\epsilon,\epsilon^{-1}[$ we put 
$$
\tau^{\epsilon}_h:=\inf\big\{t\ge 0\colon\  X^{u}_{t}
\notin[\epsilon\,,\epsilon^{-1}]\big\}\wedge h\wedge T_1. 
$$
Let $\cL^0\Phi(u):= (1/2)\sigma^2u^2 \Phi''(u)+ (au+c)\Phi'(u)$. By the It\^o formula 
 \bean
 \Phi(X^u_{\tau^{\epsilon}_h})
 &=&
 \Phi(u)+\sigma \int_{0}^{\tau^{\epsilon}_h}X^u_{s}\Phi'(X^u_{s})\,dW_s
 + \int_{0}^{\tau^{\epsilon}_h}\cL^0\Phi(X^u_{s})ds\\
 &&+
 \int_{0}^{\tau^{\epsilon}_h}\int(\Phi(X^u_{s-}+x)-\Phi(X^u_{s-}))p_P(ds,dx).
 \eean
Due to the strong Markov property $\Phi(u)=\E\,\Phi(X^u_{\tau^{\epsilon}_h})$ {\color{black}(since $\Phi(u)=0$ for $u\le 0$)}. For every $\epsilon>0$ the integrands above are 
bounded by  constants and, hence, the expectation of the stochastic integral with respect to the Wiener process is zero.  The expectation of  the integral with respect to then integer-valued measure $p_P(ds,dx)$ is equal to the integral with respect to the compensator of the latter, that is to 
$$
 \E \int_{0}^{\tau^{\epsilon}_h}\int(\Phi(X^u_{s-}+x)-\Phi(X^u_{s-}))ds\Pi_P(dx). 
$$
Moreover, $\tau^\epsilon_h=h$ when $h$ is sufficiently small  
(the threshold below which the equality holds, of course, depends on $\omega$). 

It follows that, independently of $\epsilon$,  
$$
\frac 1h \E\int_{0}^{\tau^{\epsilon}_h}\Big(\frac 12{\sigma^2}
  (X^u_{s})^2\Phi''(X^u_{s})+(aX^u_{s}+c)\Phi'(X^u_{s})\Big)ds\to \cL^0 \Phi(u) 
$$
as $h\to 0$.   Finally,  
$$
\frac 1h   \E \int_{0}^{\tau^{\epsilon}_h}\int(\Phi(X^u_{s-}+x)-\Phi(X^u_{s-}))ds\Pi_P(dx)\to 
\int (\Phi(u+y)-\Phi(u))\Pi_P(dy). 
$$ 
It follows that $\Phi$ satisfies the equation (\ref{int-diffG}). \fdem

\smallskip
\noindent 
{\sl Remark.} 
The equation  (\ref{int-diffG}) holds in the viscosity sense, i.e.  without additional assumptions on smoothness of $\Psi$, see \cite{Bel-Kab}. 


\section{Exponentially distributed jumps: from IDE to ODE}
In the case of exponentially distributed jumps 
the integro-differential equation can be written as 
\beq
\label{EIDE}
\cL\Phi(u)+
\frac {\alpha_1} {\mu_1}\int_{0} ^\infty\Phi(u-y)e^{-y/\mu_1}dy
+\frac {\alpha_2} {\mu_2}\int_{0} ^\infty\Phi(u+y)e^{-y/\mu_2}dy=0.
\eeq
where 
$$
\cL\Phi(u):= \frac 12 \sigma^2u^2 \Phi''(u)+ (au+c)\Phi'(u) -(\alpha_1+\alpha_2)\Phi(u).
$$ 
Changing variables in the integrals we get that 
\beq
\label{abr}
\cL\Phi(u)+\frac{\alpha_1} {\mu_1}I_1(u)+\frac{\alpha_2} {\mu_2}I_2(u)=0
\eeq
where 
$$
I_1(u):=\int_{-\infty} ^u\Phi(z)e^{-(u-z)/\mu_1}dz, \qquad I_2(u):=\int_{u} ^\infty\Phi(z)e^{-(z-u)/\mu_1}dz.
$$
Note that   $I'_1 =\Phi-(1/\mu_1)I_1$ and $I'_2 =-\Phi+(1/\mu_2)I_2$.
 

\smallskip
Put   
$\cT f:=\mu_1\mu_2f''+(\mu_2-\mu_1)f'-f$. 
  It  is easily seen that 
$$
\cT I_1=\mu_1\mu_2 \Phi'-\mu_1\Phi, \qquad  \cT I_2=-\mu_1\mu_2 \Phi'-\mu_2\Phi.
$$
Applying the operator $\cT$ to both sides of the equation (\ref{abr}) we get that 
 the survival probability $\Phi$ (as well as the ruin probability $\Psi:=1-\Phi$) solves the differential equation   $\cD\Phi=0$ where $\cD$ is the differential operator of the 4th order:  
\bean
\cD\Phi&:=& \cT\cL \Phi+(\alpha_1\mu_2 -\alpha_2\mu_1)\Phi'-(\alpha_1+\alpha_2)\Phi\\
&=& \mu_1\mu_2(\cL \Phi)''+(\mu_2-\mu_1)(\cL \Phi)'-\cL\Phi+(\alpha_1\mu_2 -\alpha_2\mu_1)\Phi'-(\alpha_1+\alpha_2)\Phi
\eean
After simple calculation we get that  that the obtained 4th order equation in fact is the following third order differential equation for $G:=\Phi'$:   
$$
\tilde g_3(u)G'''+\tilde g_2(u)G''+\tilde g_1(u)G'+\tilde g_0(u)G=0
$$
where the coefficients (depending on $u$) are: 
\bean
\tilde g_3(u)&:=&\frac 12 \sigma^2 \mu^2 u^2,\\ 
\tilde g_2(u)&:=& \mu^2  ((a+2\sigma^2)u+c)+\frac 12 \Delta \mu \sigma^2 u^2,\\
\tilde g_1(u)&:=&\mu^2 (2a+\sigma^2-\alpha_1-\alpha_2)+ \Delta \mu(\sigma^2u+au+c)  -\frac 12 \sigma^2 u^2,\\
\tilde g_0(u)&:=& -au-c +  \Delta \mu(a-\alpha_1-\alpha_2)+(\alpha_1\mu_2 -\alpha_2\mu_1)
\eean
with the abbreviations $\Delta \mu:=\mu_2-\mu_1$, $\mu^2:=\mu_1\mu_2$. Since  $\mu_1>0 , \mu_2>0$  we get from here the equation 
with unit coefficient at the third derivative: 
\beq
\label{G3}
G'''+q_2(u)\ G''+q_1(u)\ G'+q_0(u)\ G(u)=0
\eeq
where 
\bean
q_2(u)&:=&\frac {\Delta \mu}{\mu^2 } +\frac{2(a+2\sigma^2)}{  \sigma^2}\frac {1}{u}+\frac{2c}{  \sigma^2} \frac{1}{u^2},\\ 
q_1(u)&:=& -\frac 1{\mu^2 }+\frac{2(a+\sigma^2)\Delta \mu}{\sigma^2\mu^2 }\frac{1}{u}
+\frac{2(\Delta \mu\ c +\mu^2  (2a+\sigma^2-\alpha_1-\alpha_2))}{\mu^2  \sigma^2}\frac{1}{u^2},              \\
q_0(u)&:=& -\frac{2a}{\mu^2  \sigma^2}\frac{1}{u}+\frac{2(\Delta\mu(a-\alpha_1-\alpha_2)+(\alpha_1\mu_2 -\alpha_2\mu_1)-c)}{\mu^2  \sigma^2}\frac{1}{u^2}.  
\eean

Let denote $\cJ$ the operator in the left-hand side of the   basic IDE acting in the space of sufficiently smooth functions.  
Then $\cT\cJ=\cD$.  Note that  ${\rm dim}\,{\rm Ker}\, \cD=4$. Also, we have the inclusion 
${\rm Ker}\, \cJ\subseteq {\rm Ker}\, \cD$. 
The kernel of the 2nd order differential operator $\cT$ is a 2-dimensional linear subspace generated by the 
functions $h_1(u):=e^{-u/\mu_1}$, $h_2(u):=e^{u/\mu_2}$. 
Let $f_1$ and $f_2$ be solutions of the equations $\cJ f_j=h_j$. Then 
$f_1$ and $f_2$ are linearly independent and belong to  ${\rm Ker}\, \cD$. 
The 4 functions: the identity, the survival probability $\Phi$ (assumed not to be a constant), $f_1$, $f_2$ form a basis in ${\rm Ker}\, \cD$. Indeed,  if  their linear combination 
$$
a_1f_1+a_2f_2+a_3 1+a_4\Phi=0,
$$
then $a_1\cJ f_1+a_2\cJ f_2=0$. That is, $a_1 h_1+a_2 h_2=0$ and $a_1=a_2=0$. But the equality  $a_3+a_4\Phi=0$ holds only if $a_3=a_4=0$.

\section{Asymptotic analysis of the differential equation for the survival probability}

We analyze the  behavior of solutions of the equation (\ref{G3})  using a result on systems with asymptotically constant coefficients. To this end, we put $y=(y_1,y_2,y_3):=(G,G',G'')$. Using the matrix notation where the vectors are columns we get from  (\ref{G3}) that $y'=A(u)y$ where 
$$
A(u):=
\begin{bmatrix}
0 & 1 &0\\ 
0 & 0 &1\\
- q_0(u) & - q_1(u) &- q_2(u)  
\end{bmatrix}
, \qquad 
A(\infty):=
\begin{bmatrix}
0 & 1 &0\\ 
0 & 0 &1\\
0 & 1/\mu^2 &- (\Delta \mu)/\mu^2  
\end{bmatrix}.
$$
Then $A(u)=A(\infty)+V(u)$ where the matrix $V(u):=A(u)-A(\infty)$ is such that its  norm  
$|V'(u)|$ (the Euclidean or any other) is integrable on $[1,\infty[$.

Let  $\lambda_j=\lambda_j(u)$, $j=1,2,3$, be the roots  of the characteristic equation
\beq
\label{cubic}
\lambda^3+q_2(u)\ \lambda^2+q_1(u)\ \lambda+q_0(u)=0.
\eeq
Note that 
\beq
\label{3.}
 \lambda_1+\lambda_2+\lambda_3  =  -q_2(u),\quad {\lambda_1\lambda_2}+{\lambda_2 \lambda_3}+{\lambda_1\lambda_3}  =  q_1(u), \quad \lambda_1\lambda_2 \lambda_3  =  -q_0(u).
\eeq


Recall  that we are working under the assumptions  $a>\sigma^2/2>0$ and $\mu_1,\mu_2>0$.  
Then $q_0(u)\to 0$,  
$$
q_2(u) \to \frac {\Delta \mu}{\mu^2 }\neq 0, \quad q_1(u) \to  -\frac {1}{\mu^2 }<0,  \quad u q_0(u) \to  -\frac{2a}{\mu^2  \sigma^2}<0, \quad u\to \infty.  
$$
The Cardano formulae imply that the roots $\lambda_j(u)$ are continuous functions  having finite limits as $u\to \infty$.  
According to the last equation in (\ref{3.}), at least one root, say, $\lambda_3(u)$ tends to zero as $u\to \infty$. Since $q_1(\infty)\neq 0$, two  other roots have nonzero  limits satisfying   the system 
$$ \lambda_1(\infty)+\lambda_2(\infty)=-q_2(\infty), \qquad
\lambda_1(\infty)\lambda_2(\infty)=q_1(\infty),  
$$
that is, $\lambda_1(\infty)=1/\mu_2$, $\lambda_2(\infty)=-1/\mu_1$. Thus,  we obtain that 
\beq
\lambda_1(u)=\frac 1{\mu_2}+o(1), \quad \lambda_2(u)=- \frac 1{\mu_1}+o(1), \quad  \lambda_3(u)=-\frac {2a}{\sigma^2}\frac 1u+o(u^{-1}). 
\eeq 
Applying  the implicit function theorem we conclude that the function $\lambda_3(u)$ can be expanded in powers of $u^{-1}$ In particular,  
\beq
\label{O}
 \lambda_3(u)=-\frac {2a}{\sigma^2}\frac 1u+O(u^{-2}), \quad u\to \infty.
\eeq 
The numbers $\lambda_1(\infty)=1/\mu_2$, $\lambda_2(\infty)=-1/\mu_1$, and $\lambda_3(\infty)=0$ are eigenvalues of the matrix $A(\infty)$. 

The conditions of Th. VII-5-3 from \cite{HS} are fulfilled and, therefore, the fundamental matrix of the equation $y'=A(u)y$ has the form 
$$
P_0(u)(I+H(u))\exp\left\{\int ^u \Lambda (s)ds \right\}
$$
where the  matrix-valued functions $P_0(u)$ and $H(u)$ are continuous, $P_0(u)\to Q$, 
$H(u)\to 0$ as $u\to \infty$, $\Lambda (s):={\rm diag}\,(\lambda_1(s),\lambda_2(s),\lambda_3(s))$,  the  columns of $Q$ are eigenvectors of $A(\infty)$ corresponding to the eigenvalues $\lambda_1(\infty)$, $\lambda_2(\infty)$, $\lambda_3(\infty)$,  that is 
$$
Q=
\begin{bmatrix}
1 & 1 & 1\\ 
1/\mu_2 & -1/\mu_1 &0\\
1/\mu_2^2 & 1/\mu_1^2 &0  
\end{bmatrix}.
$$ 

A general solution of solution $G$ of (\ref{G3}) is a linear combination of functions 
\bean
h_1(u)&:=&(1+\theta_1(u)) \exp\left \{\int_1^u \left (\frac 1{\mu_2}+\gamma_1(s)\right)ds\right\}, \\ 
h_2(u)&:=&(1+\theta_2(u)) \exp\left \{\int_1^u \left (-\frac 1{\mu_1}+\gamma_2(s)\right)ds\right\}, \\
h_3(u)&:=&(1+\theta_3(u)) \exp\left \{-\frac {2a}{\sigma^2}\ln u\right\} \exp\left \{\int_1^u \gamma_3(s) ds\right\},
\eean
where continuous functions $\theta_j$ and $\gamma_j$ vanish at infinity for $j=1,2,3$ and 
the function $\gamma_3$ is integrable. 
 Thus, the general solution of the 4th order differential  equation is a linear combination of a constant function and  these three functions above.  In particular, the ruin probability $\Psi$ is given by a certain linear combination. Clearly, the latter cannot involve the unbounded function corresponding to the integral of first function. The integrals of others are  bounded and, therefore, 
$$
\Psi(u)= c_0+ c_2H_2(u)+c_3H_3(u)
$$  
where 
$$
H_j(u):=\int_u^\infty h_j(s)ds, \quad j=2,3.  
$$
Note that $c_0=0$, the  integral $H_2(u)$  converges to zero exponentially fast, 
and 
$$
\frac {H_3(u)}{u^{1-2a/\sigma^2}}\to \frac 1{2a/\sigma^2-1}\exp\left \{\int_1^\infty \gamma_3(s) ds\right\}\neq 0
$$
 as $u\to \infty$.  
in virtue of Proposition \ref{lab} $c_3\neq 0$ and we easily obtain that $\Psi(u)\sim C u^{-\beta}$ where $C>0$ and $\beta=2a/\sigma^2-1$.      

\smallskip
{\bf Acknowledgement.} This work was supported by the Russian Science Foundation grant  
20-68-47030 and was completed during the stay of the  first author  at FRIAS.


 \end{document}